\documentclass[12pt,english]{article}
\usepackage[margin=0.75in]{geometry}
\usepackage[utf8]{inputenc}
\usepackage[english]{babel}
\usepackage{amsthm}
\usepackage{graphicx}
\usepackage{url}
\usepackage{hyperref}
\usepackage{amsmath}
\usepackage{listings}
\usepackage{amsfonts}
\usepackage{xcolor}
\usepackage{float}
\newtheorem{theorem}{Theorem}

\def\Xint#1{\mathchoice
   {\XXint\displaystyle\textstyle{#1}}%
   {\XXint\textstyle\scriptstyle{#1}}%
   {\XXint\scriptstyle\scriptscriptstyle{#1}}%
   {\XXint\scriptscriptstyle\scriptscriptstyle{#1}}%
   \!\int}
\def\XXint#1#2#3{{\setbox0=\hbox{$#1{#2#3}{\int}$}
     \vcenter{\hbox{$#2#3$}}\kern-.5\wd0}}
\def\dashint{\Xint-}

\begin{document}

\title{Analytical continuation of prime zeta function for $\Re(s)>\frac{1}{2}$ assuming (RH)}

\author{Artur Kawalec}

\date{}
\maketitle

\begin{abstract}
We derive a simple expression to analytically continue the prime zeta function to the domain $\Re(s)>\frac{1}{2}$ assuming (RH) and taking into account a proper branch cut. We also verify the formula numerically and provide several plots.
\end{abstract}

\section{Introduction}
Let $p=\{2,3,5,7,\ldots\}$ be a sequence of primes and the prime zeta function

\begin{equation}\label{eq:1}
P(s)=\sum_{p}\frac{1}{p^s}
\end{equation}
is defined as a generalized zeta series over primes, which converges absolutely for $\Re(s)>1$. Another representation of the prime zeta function is the very well-known M\"obius inversion formula

\begin{equation}\label{eq:1}
P(s)=\sum_{n=1}^{\infty}\frac{\mu(n)}{n}\log\zeta(ns)
\end{equation}
which further continues $P(s)$ for $\Re(s)>0$ (see [2]). And also, there is a natural boundary line at $\Re(s)=0$, after which $P(s)$ cannot be further continued to $\Re(s)\leq 0$ because of a dense accumulation of singularities that form near $\Re(s)=0$ [5, p. 215]. In our previous article [3] we derived a series expansion of the prime zeta about $s=1$ using Stieltjes integration. In this article, we also show by similar Stieltjes integration, that

\begin{theorem}
\label{thm:fta}

\begin{equation}\label{eq:1}
P(s)=\lim_{x\to \infty}\Bigg\{\sum_{p\leq x}\frac{1}{p^s}+\operatorname{E_1}[(s-1)\log(x)]\Bigg\}
\end{equation}
which is analytic for $\Re(s)>\tfrac{1}{2}$ with a branch cut on $(\tfrac{1}{2},1]$ assuming (RH).

\end{theorem}

\noindent
And the exponential integral is defined by

\begin{equation}\label{eq:1}
\operatorname{E_1}(z)=\int_{z}^{\infty}\frac{e^{-t}}{t}dt, \quad |\arg(z)|<\pi
\end{equation}
valid for complex $z$ with a branch cut on $(-\infty, 0]$, and its series expansion is

\begin{equation}\label{eq:1}
\operatorname{E_1}(z)=-\gamma-\log(z)-\sum_{k=1}^{\infty}\frac{(-z)^k}{k\cdot k!}
\end{equation}

\noindent
We now prove this Theorem in the next Section.

\section{Proof of Theorem 1}

To show that, let us review some basic definitions. The prime counting function is defined by

\begin{equation}\label{eq:1}
\pi(x)=\sum_{p\leq x} 1
\end{equation}
for positive integer argument $x>0$, and by the average value

\begin{equation}\label{eq:1}
\pi(x)=\frac{1}{2}\Bigg[\sum_{p<x} 1 +\sum_{p\leq x} 1\Bigg]
\end{equation}
for positive real argument $x>0$. This means that its value is defined as an average of the two sides of the step when $x$ is an integer. The asymptotic formula is

\begin{equation}\label{eq:1}
\pi(x)=\operatorname{li}(x)+f(x)
\end{equation}
where the logarithmic integral is defined by the P.V. integral function

\begin{equation}\label{eq:1}
\operatorname{li}(x)=\dashint_{0}^{x}\frac{1}{\log t} dt
\end{equation}
for real $x>0$. One also has the relation
\begin{equation}\label{eq:1}
\frac{d}{dx}\operatorname{li}(x)=\frac{1}{\log x}.
\end{equation}
And as for the remainder $f(x)$ the best estimate is
\begin{equation}\label{eq:1}
f(x)=O(\sqrt{x}\log x)
\end{equation}
assuming (RH).
We prove Theorem $1$ by Stieltjes integration of the remainder term as follows:

 \begin{equation}\label{eq:1}
\begin{aligned}
P(s)
&= \sum_{p\leq x}\frac{1}{p^s}+\int_{x}^{\infty}\frac{1}{t^{s}}\, d\pi(t) \\[1.2em]
& = \sum_{p\leq x}\frac{1}{p^s}+\int_{x}^{\infty} \frac{1}{t^{s}} d\operatorname{li}(x)\,  + \int_{x}^{\infty} \frac{1}{t^{s}} df(t) \\[1.2em]
&= \sum_{p\leq x}\frac{1}{p^s}+\int_{x}^{\infty} \frac{1}{t^{s}\log(t)} dt\,  + \Biggl[\frac{f(t)}{t^s}\Biggr]_{x}^\infty+ s\int_{x}^\infty t^{-s-1}f(t)\, dt \\[1.2em]
&=\sum_{p\leq x}\frac{1}{p^s}+\operatorname{E_1}[(s-1)\log(x)]-\frac{f(x)}{x^s}+s\int_{x}^\infty t^{-s-1}f(t)\, dt \\[1.2em]
&=\sum_{p\leq x}\frac{1}{p^s}+\operatorname{E_1}[(s-1)\log(x)]-\frac{\pi(x)-\operatorname{li}(x)}{x^s}+O\Big(x^{\tfrac{1}{2}-s}\log x\Big)\,  \\[1.2em]
\end{aligned}
\end{equation}

\noindent
The integral can be expressed in terms of the exponential integral as

\begin{equation}\label{eq:1}
\int_{x}^{\infty} \frac{1}{t^{s}\log t} dt=\operatorname{E_1}[(s-1)\log(x)]
\end{equation}
by substituting $u=\log t$ and $t=e^u$ and $dt=e^u du$. And also the size of the boundary term

\begin{equation}\label{eq:1}
O\Big(\frac{\pi(x)-\operatorname{li}(x)}{x^s}\Big) = O\Big(x^{\tfrac{1}{2}-s}\log x\Big)
\end{equation}
is also the same as the last integral in (12). As a result, when taking the limit $x\to \infty$, we have that

\begin{equation}\label{eq:1}
P(s)=\lim_{x\to \infty}\Bigg\{\sum_{p\leq x}\frac{1}{p^s}+\operatorname{E_1}[(s-1)\log(x)]-\frac{\pi(x)-\operatorname{li}(x)}{x^s}\Bigg\}
\end{equation}
for $\Re(s)>\tfrac{1}{2}$, but the boundary term (14) doesn't improve convergence of this equation, hence it may be dropped. And we must introduce at branch cut $(\tfrac{1}{2},1]$  because of log in the exponential integral (see series expansion (5)).

\section{Numerical Computation}
In this Section, we write a simple script in Pari/GP [4] to compute the prime zeta by equation (3), where we utilize the built-in exponential integral function $E_1(z)$ as $\textbf{eint1(z)}$. In Fig. 1 we plot the prime zeta equation (2) and compare with equation (3) for real $s$ domain and limit variable $x=10^4$,  but noting that we must compute $\Re[ P(s)]$ since it's on the branch cut. We also see that there is a near perfect reproduction, but near $s=0.5$ the convergence of equation (3) starts deviating. And in Fig. 2-3, we plot the same functions but on the vertical line with real part $\sigma=0.75$ and $t=0.1$ to $50$ for real and imaginary part of $P(s)$ with limit variable $x=10^4$. We also see a perfect match.

\lstset{language=C,deletekeywords={for,double,return},caption={A sample Pari/GP code for calculating prime zeta and plotting Fig. 1-3},label=DescriptiveLabel,captionpos=b,showstringspaces=false}
\begin{lstlisting}[frame=single]
\\ Define P(s) by Equation (2) valid for Re(s)>0
PrimeZeta(s)=sum(n=1, 10^3, moebius(n)/n*log(zeta(n*s)));

x = 10^3; \\ Set value for limit variable
a = 0.75; \\ set real part of vertical line

\\ Define P(s) by Equation (3) valid for Re(s)>1/2
PrimeZeta_Eq3(s)=sum(n=1,primepi(x),1.0/prime(n)^s)+eint1((s-1)*log(x));

\\ Fig. 1 Plot
ploth(s=0.5001, 2, [real(PrimeZeta(s)),real(PrimeZeta_Eq3(s))])

\\ Fig. 2 real plot
ploth(t=0.1, 50, [real(PrimeZeta(a+I*t)),real(PrimeZeta_Eq3(a+I*t))])

\\ Fig. 3 imag plot
ploth(t=0.1, 50, [imag(PrimeZeta(a+I*t)),imag(PrimeZeta_Eq3(a+I*t))])
\end{lstlisting}

\texttt{Email: art.kawalec@gmail.com}

\begin{figure}[H]
\centering
  \renewcommand{\thefigure}{1}%
  \includegraphics[width=160mm]{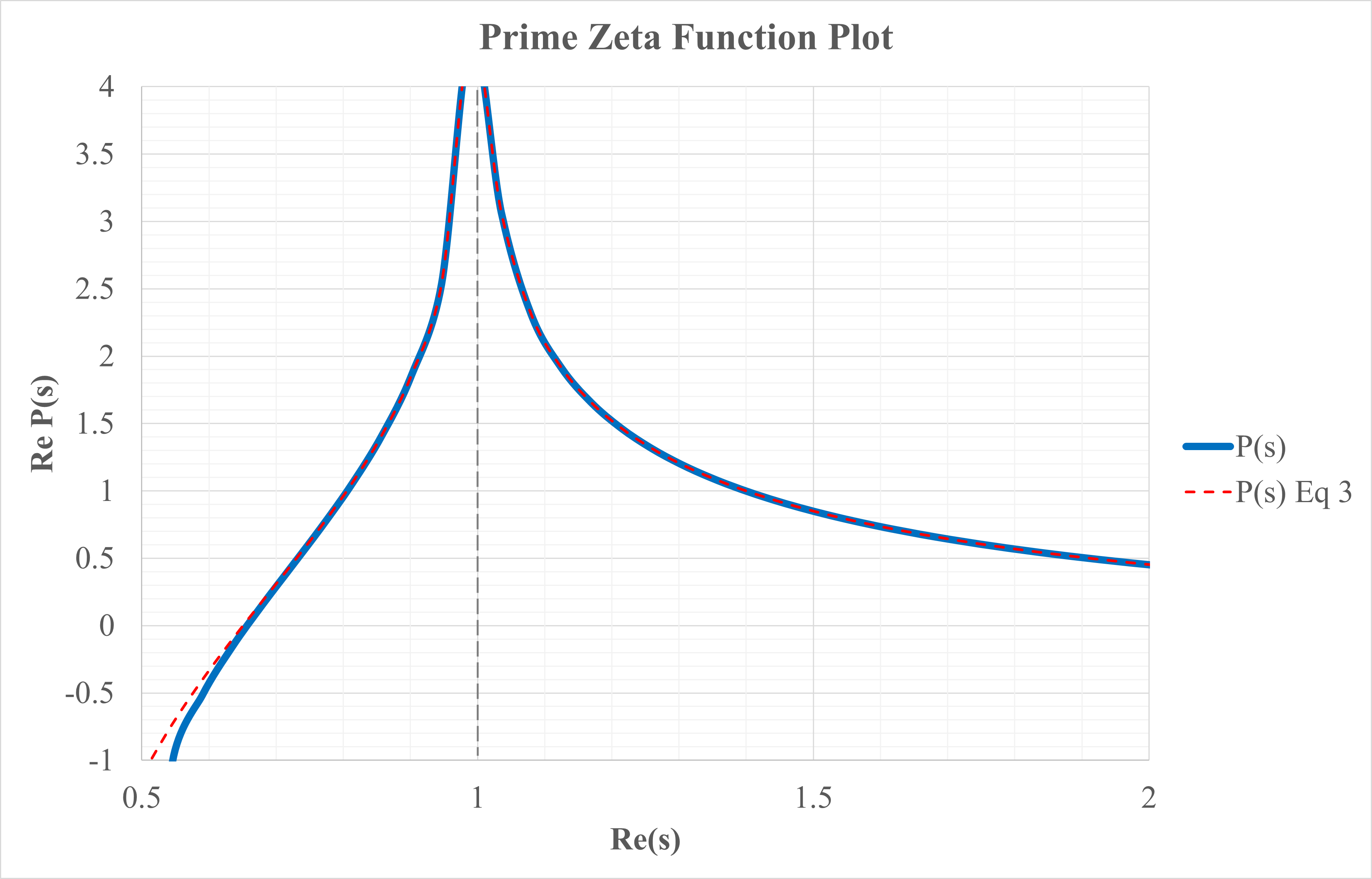}\\
  \caption{A plot of $\Re[P(s)]$ by equation (3) for $\Re(s)$ variable for limit variable $x=10^4$}\label{1}
\end{figure}

\begin{figure}[H]
  \centering
  \renewcommand{\thefigure}{2}%
  \includegraphics[width=160mm]{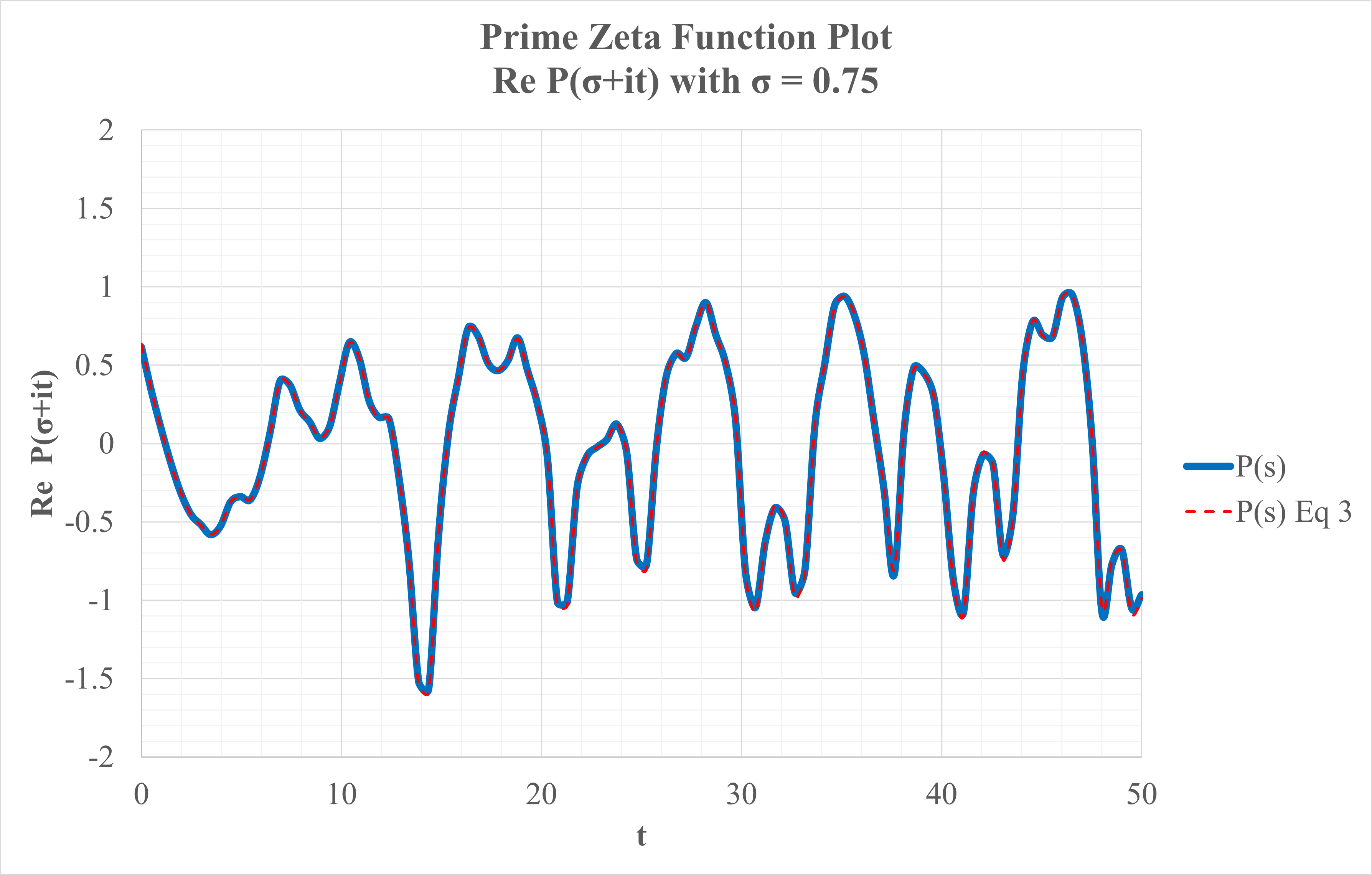}\\
  \caption{A plot of $\Re[P(\sigma+it)]$ by equation (3) for at $\sigma=0.75$ and vertical line $t=0.1$ to $t=50$ for limit variable $x=10^4$}\label{1}
\end{figure}

\begin{figure}[H]
   \centering
   \renewcommand{\thefigure}{3}%
  \includegraphics[width=160mm]{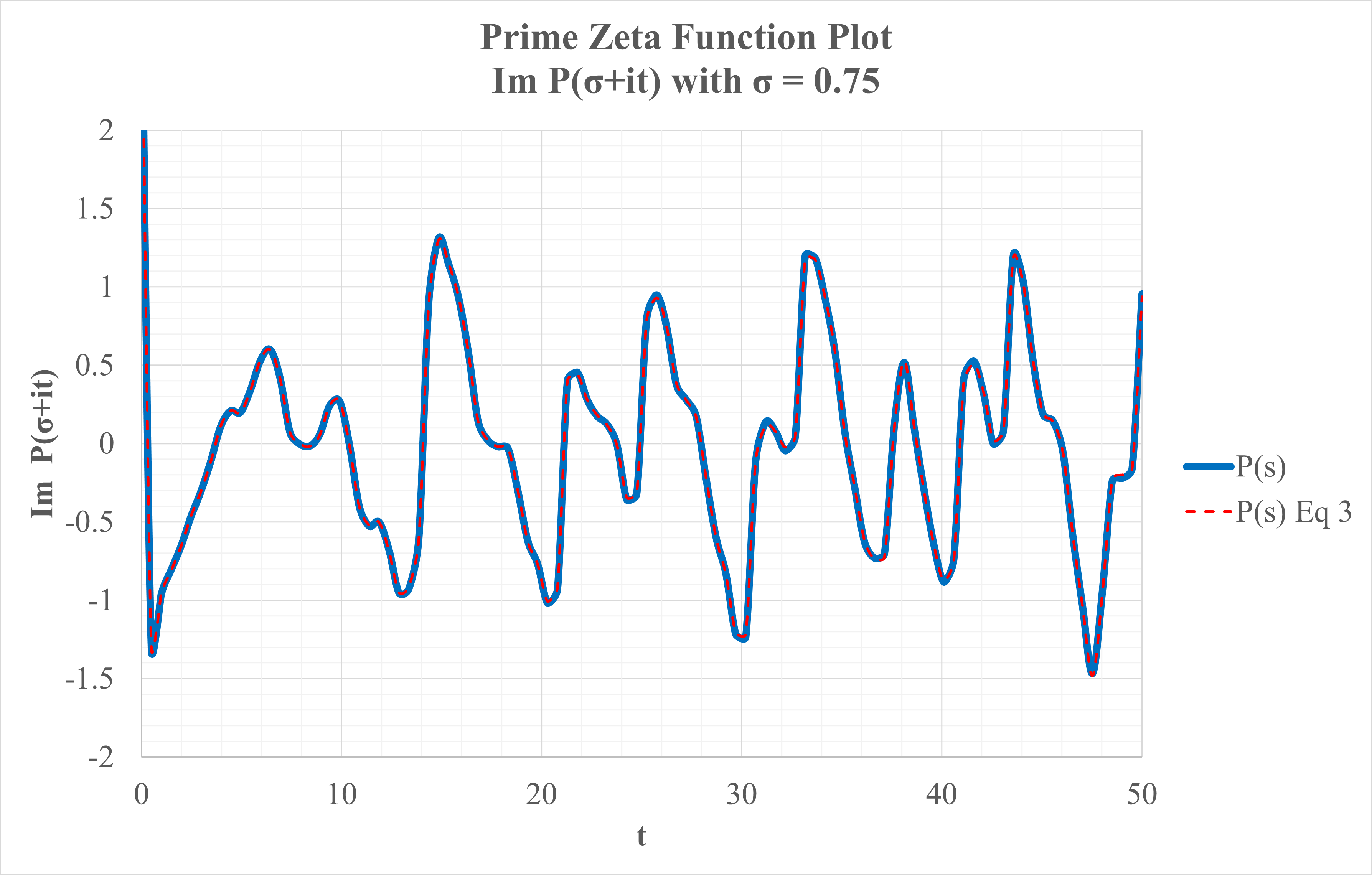}\\
  \caption{A plot of $\Im[P(\sigma+it)]$ by equation (3) for at $\sigma=0.75$ and vertical line $t=0.1$ to $t=50$ for limit variable $x=10^4$}\label{1}
\end{figure}

\end{document}